\documentclass[11pt]{article}
\usepackage{amsfonts}
\newcommand{\beq}{\begin{equation}}
\newcommand{\eeq}{\end{equation}}
\newcommand{\beqa}{\begin{eqnarray}}
\newcommand{\eeqa}{\end{eqnarray}}
\newcommand{\beaa}{\begin{eqnarray*}}
\newcommand{\eaa}{\end{eqnarray*}}
\newcommand{\nn}{\hfill\nonumber}
\newtheorem{thm}{Theorem}
\newtheorem{lem}{Lemma}
\newtheorem{prop}{Proposition}

\newcommand \nc {\newcommand}
\nc \bde[1]{ \medskip \noindent {\bf{Definition #1}} }
\nc \bex[1]{ \medskip \noindent {\it{Example #1}} }
\nc \bre[1]{ \medskip \noindent {\it{Remark #1}} }
\nc \ede { \medskip }
\nc \eex { \medskip }
\nc \ere { \medskip }
\nc \eqref[1] {(\ref{#1})}
\nc \thref[1]{Theorem \ref{t#1}}   
\nc \leref[1]{Lemma \ref{l#1}}
\nc \prref[1]{Proposition \ref{p#1}}
\nc \coref[1]{Corollary \ref{c#1}}
\nc \deref[1]{Definition #1}
\nc \exref[1]{Example #1}
\nc \reref[1]{Remark #1}
\newcommand {\normprod}[1]{ {\textrm{:}}{#1}{\textrm{:}} } 
\def \W {W_{1+\infty}}
\def \WN {\W(N)}
\def \a {\alpha}

\def \A {{\mathcal A}}
\def \D {{\mathcal D}}
\def \M {{\mathcal M}}
\def \Cset {{\mathbb C}}
\def \Zset {{\mathbb Z}}
\def \Nset {{\mathbb N}}
\def \rank { {\mathrm{rank}} }
\def \mult { {\mathrm{mult}} }
\def \span { {\mathrm{span}} }
\def \const { {\mathrm{const}} }
\renewcommand \ker { {\mathrm{Ker}} }
\begin{document}
\title{{\LARGE\bf{Highest weight modules of $W_{1+\infty}$, Darboux
transformations and the bispectral problem
\thanks{Presented at the conference Geometry and Mathematical Physics,
Zlatograd 95.}
}}}
\author{
B.~Bakalov
\thanks{E-mail: bbakalov@fmi.uni-sofia.bg}
\quad
E.~Horozov
\thanks{E-mail: horozov@fmi.uni-sofia.bg}
\quad
M.~Yakimov
\thanks{E-mail: myakimov@fmi.uni-sofia.bg}
\\ \hfill\\ \normalsize \textit{
Department of Mathematics and Informatics, }\\
\normalsize \textit{Sofia University, 5 J. Bourchier Blvd.,
Sofia 1126, Bulgaria }     }
\date{}
\maketitle
\vspace{-8cm}
\begin{flushright}
{\tt{ q-alg/9601017 }}
\end{flushright}
\vspace{7cm}
\setcounter{section}{-1}
\section{Introduction}
In this paper we announce the existence of a large class of {\em{bispectral\/}}
ordinary differential operators and a series of their properties. Following
\cite{DG} we call an operator $L(x,\partial_x)$ bispectral if it possesses a
family of eigenfunctions $\Psi(x,z),$ which are also eigenfunctions for another
differential operator $\Lambda(z,\partial_z),$ but this time in the ``spectral
parameter'' $z$, to wit
    \beqa
&&L(x,\partial_x) \Psi(x,z) = f(z) \Psi(x,z) \label{0.1} \hfill\\
&&\Lambda(z,\partial_z) \Psi(x,z) = \Theta(x) \Psi(x,z) \label{0.2} \hfill
     \eeqa
for some functions $f(z),\Theta(x).$

The problem of describing bispectral operators has its roots in several
mathematical and physical issues, e.g.\ -- computer tomography. For more
motivation and background we recommend \cite{DG, G1, G3}.

The first general result in the direction of classifying bispectral operators
belongs to J.~J.~Duistermaat and F.~A.~Gr\"unbaum \cite{DG} who  determined
all second order bispectral operators $L$.
Their answer is as follows. If we write the operator $L$ as
  $$
L = \partial_x^2 + u(x),
  $$
the bispectral potentials $u(x)$ are given (up to translations and rescalings
of $x$ and $z$) apart from the obvious Airy ($u(x)=x$) and Bessel
($u(x)=cx^{-2}$) potentials  by potentials obtained by finitely many
``rational'' Darboux transformations from $u(x) = 0$  and
$u(x) = (-1/4)x^{-2}$.

Every operator $L(x,\partial_x)$ can be considered as an element of a maximal
algebra $\A$ of commuting ordinary differential operators \cite{BC}. Led by
this observation G.~Wilson \cite{W} introduced the following terminology.
He called such an algebra {\it {bispectral}} if there exists a
joint eigenfunction $\Psi(x,z)$ for the operators $L$ in $\A$ that satisfy also
the eq. \eqref{0.2}. The dimension of the space of eigenfunctions $\Psi(x,z)$
is called a {\it {rank}} of the commutative algebra $\A$ (see e.g.\
\cite{KrN}).
This number coincides with the greatest common divisor of the orders of the
operators in $\A$. Wilson classified all rank 1 bispectral algebras
\cite{W}.

Our construction puts in a general context all previously known results
and explains them (to some extent) from a representation-theoretic point of
view.
The main results of the present paper are as follows.

First, starting with certain generalizations of Bessel functions we construct
highest weight modules $\M_\beta$ of the Lie algebra $\W$ with highest weight
vectors -- tau-functions.

Second, performing a large class of Darboux transformations on polynomials of
the corresponding to these tau-functions differential operators (Bessel
operators) we construct families of bispectral algebras of any rank $N$.
There is a similar construction for the higher order generalizations of Airy
operators.

Third, we show that the manifolds of bispectral operators obtained by
Darboux transformations on powers of Bessel operators are in one to one
correspondence with the manifolds of tau-functions lying in  the modules
$\M_\beta$. An immediate corollary is that they are preserved by hierarchies
of symmetries generated by subalgebras of $\W$.

At the end we point out that the suggested method allows to obtain the
bispectral algebras algorithmically despite that we use highly transcendental
functions like Bessel and Airy ones. We covered all classes and examples of
bispectral operators which we know from the literature
\cite{DG, W, Z, G3}, etc.
We conjecture that all bispectral scalar differential operators can be obtained
in this way.
\\ \hfill \\
The proofs of the results presented here are performed in
\cite{BHY}--\cite{BHY3}.
\section{ Darboux transformations and bispectral algebras}
The framework of our construction is Sato's theory of KP--hierarchy \cite{S,
SW, KRa}. In particular our eigenfunctions are Baker (or wave) functions
$\Psi_V(x,z)$ corresponding to points (or planes) $V$ in Sato's Grassmannian
$Gr$. We obtain our bispectral algebras by applying a version of {\it {
Darboux transformations}} on specific wave functions
which we call {\it {Bessel {\rm{(and}} Airy\/{\rm{)}} wave functions}}.

For $\beta \in \Cset^N$ such that
 $\sum_{i=1}^{N}\beta_i = N(N-1)/2$
we introduce the ordinary differential operator
\beq
  P_{\beta}(D_z) = (D_z-\beta_1)(D_z-\beta_2) \cdots (D_z-\beta_N)
  \label{2.13}
\eeq
where $D_z = z\partial_z$ and consider the differential equation
\beq
 P_{\beta} (D_z) \Phi_\beta (z) = z^N \Phi_\beta (z).
 \label{2.14}
\eeq
For every sector $S$ with a center at the irregular singular point $z= \infty$
and an angle less than $2 \pi$ the equation \eqref{2.14} has a solution
$\Phi_\beta$
with an asymptotics
\beq
 \Phi_\beta(z) \sim\Psi_{\beta}(z) =
e^z \Bigl( 1+ \sum_{k = 1}^\infty a_k (\beta) z^{-k} \Bigr)
\label{2.14'5}
\eeq
for $|z| \to \infty$, $z \in S$. Here
$a_k(\beta)$ are symmetric polynomials in $\beta_i$.
The function $\Phi_\beta(z)$ can be taken to be (up to a rescaling) the Meijer's
 $G$-function
$G^{N0}_{0N}\bigl((-z/N)^N \big| (1/N)\beta\bigr)$
(see \cite{BE}, \S $5.3$).
The next definition is fundamental for the present paper.

\bde{1}
{\it {Bessel wave function}} is called the function
$ \Psi_{\beta}(x,z) = \Psi_{\beta}(xz) $ (cf. \cite{F, HH, Z}).
The {\it {Bessel operator}} $L_{\beta}$ is defined as
\beq
 L_{\beta}(x,\partial_x) = x^{-N} P_{\beta} (D_x).
 \label{2.16}
\eeq
A Bessel wave function $\Psi_\beta$ defines a plane $V_\beta \in Gr$ (called
{\em{Bessel plane\/}}) by the standard procedure:
$$ V_\beta = \span \{ \partial_x^k \Psi_\beta(x,z) |_{x = 1}\}. $$
We denote by $\tau_\beta(t)$ the corresponding tau-function and call it
{\em{Bessel tau--function\/}}.
\ede

Classically, a Darboux transformation \cite{Da} of a differential operator $L$
presented as a product $L = Q P$ is defined by exchanging the places of the
factors, i.e.\ $ \overline L = P Q $. Obviously, if $\Psi(x,\lambda)$ is an
eigenfunction of $L$, i.e.\
$L(x,\partial_x)\Psi(x,\lambda) = \lambda\Psi(x,\lambda)$
then $P\Psi(x,\lambda)$ is an eigenfunction of $\overline L$.
Our definition of a Darboux transformation puts the emphasis rather on the
eigenfunctions $\Psi(x, \lambda)$ and $P \Psi(x, \lambda)$ than on the
operators $L$ and $\overline L$.

\bde{2}
We say that a plane $W$ (or the corresponding wave
function $\Psi_W(x,z)$) is a {\em Darboux transformation\/} of the plane $V$
(respectively wave function $\Psi_V(x,z)$) iff there exist monic polynomials
$f(z)$, $g(z)$ and differential operators $P(x,\partial_x)$, $Q(x,\partial_x)$
such that
\beqa
&&\Psi_W(x,z)=\frac{1}{g(z)} P(x,\partial_x) \Psi_V(x,z),
\label{3.8} \\
&&\Psi_V(x,z)=\frac{1}{f(z)} Q(x,\partial_x) \Psi_W(x,z).
\label{3.9}
\eeqa
\ede
Simple consequences of \deref{2} are the identities
\beqa
&&PQ\Psi_W(x,z)=f(z)g(z)\Psi_W(x,z),
\label{3.10}\\
&&QP\Psi_V(x,z)=f(z)g(z)\Psi_V(x,z).
\label{3.11}
\eeqa
The operator $\overline L=PQ$ is a Darboux transformation of $L=QP$.
Obviously \eqref{3.8} implies the inclusion
\beq
gW\subset V.
\label{3.12}
\eeq
Conversely, if \eqref{3.12} holds there exists $P$ satisfying \eqref{3.8}.
Therefore $W$ is a Darboux transformation of $V$ iff
\beq
fV\subset W\subset \frac1gV
\label{3.13}
\eeq
for some polynomials $f(z)$, $g(z)$.

Obviously some of the Bessel planes $V_\beta, \; \beta \in \Cset^N$ can be
obtained by a Darboux transformation from other Bessel planes $V_{\beta'}, \;
\beta' \in \Cset^{N'}, \;
N'<N.$ Of course, we are interested in $\beta$ which do not have this property.
We call them {\em{generic}}.

To each plane $W$ one can associate the {\em spectral algebra\/} $A_W$ of
polynomials $f(z)$ that leave $W$ invariant.
For each $f(z)\in A_W$ one can show that there exists a unique
differential operator $L_f(x,\partial_x)$, the order of $L_f$ being equal to
the degree of $f$, such that
\beq
L_f(x,\partial_x) \Psi_W(x,z) = f(x) \Psi_W(x,z).
\label{1.20'5}
\eeq
We denote the commutative algebra of these operators by $\A_W$.
\begin{prop}\label{p5.4}
 For a generic $\beta\in \Cset^N$ we have
\beq
A_{V_\beta} = \Cset[z^N], \quad \A_{V_\beta} = \Cset[L_\beta].
\label{5.10}
\eeq
\end{prop}
  For a Darboux transformation of a Bessel plane $V=V_\beta$
with generic $\beta\in\Cset^N$ \prref{5.4} and \eqref{3.11} imply
\beqa
&&f(z)g(z)=h(z^N),
\hfill \label{5.14} \\
&&Q P=h(L_\beta)
\hfill \label{5.15}
\eeqa
for some polynomial $h$. The operator $P$ is determined by its kernel which is
a subspace of $\ker h(L_\beta)$. The latter is described in the following lemma
(see e.g.\ \cite{I}).
\begin{lem}\label{l5.1}
 Let $h(z)$ be a polynomial
\beq
h(z)=z^{d_0}\left(z-\lambda_1^N\right)^{d_1}\cdots
\left(z-\lambda_r^N\right)^{d_r},\quad
\lambda_i^N\ne \lambda_j^N,\ \lambda_0=0,\ d_i\ge0.
\label{5.1}
\eeq
Then we have

{\em(i)} $\ker h(L_\beta)=\bigoplus_{i=0}^r \ker
\left(L_\beta-\lambda_i^N\right)^{d_i}$.

{\em(ii)} $(L_\beta)^d=L_{\beta^d}$,  where
\beq
\beta^d=(\beta_1,\beta_1+N,\ldots,\beta_1+(d-1)N,\ldots,\beta_N,\ldots,\beta_N
+(d-1)N).
\label{5.1'5}
\eeq

{\em(iii)}  If
$\{\beta_1,\ldots,\beta_N\}=\{\underbrace{\alpha_1,\ldots,\alpha_1}_{k_1}
,\ldots,
\underbrace{\alpha_s,\ldots,\alpha_s}_{k_s} \}$
with distinct $\alpha_1, \ldots,\alpha_s$, then
$$
\ker L_\beta=\span\left\{ x^{\alpha_i}(\ln x)^k\right\}_{1\le i\le s,\ 0\le
k\le k_i-1}.
$$

{\em(iv)}  For $\lambda\ne 0$
$$
\ker\left(L_\beta-\lambda^N\right)^d =\span\left\{
\partial_\lambda^k\Psi_\beta(x,\lambda\varepsilon^j)\right\}_{0\le k\le d-1,\
0\le j\le N-1},
$$
where $\varepsilon=e^{2\pi i/N}$ is an $N$-th root of unity.
\end{lem}
We call $\ker P$ {\em{homogeneous\/}} and $\Zset_N$-{\em{invariant\/}} iff it
has a basis which is a union of:

(i) Several groups of elements supported at $0$ of the form:
\beq
\partial_y^l \Big(\sum_{k=0}^{k_0}\sum_{j=0}^{\mult(\beta_i+kN)-1}
b_{kj} x^{\beta_i + kN} y^j \Big) \Big|_{y = \ln x}, \quad 0 \leq l \leq j_0,
\label{5.20a}
\eeq
where $\mult(\beta_i + kN):=$ multiplicity of
$\beta_i+ kN$ in $\bigcup_{j=1}^N \{\beta_j+ N\Zset_{\geq0}\}$ and
$j_0 =\max\{j | b_{kj}\not=0 {\textrm{ for some }} k\}$.

(ii) Several groups of elements supported at the points $\varepsilon^i \lambda$
($0 \le i \le N-1$, $\lambda\not=0$) of the form:
\beq
\sum_{k=0}^{k_0} a_k \varepsilon^{ki}
\partial_z^k\Psi_\beta(x,z)|_{z=\varepsilon^i\lambda},\quad 0\le i\le N-1.
\label{5.20b}
\eeq
Instead of \eqref{5.20b} we can also take
\beq
\hspace*{0.5cm}
\sum_{k=0}^{k_0} a_k D^k_z\Psi_\beta(x,z)|_{z=\varepsilon^i\lambda},
 \quad 0\le i\le N-1.
\label{5.20c}
\eeq
Denote by $n_0$ the number of elements of the form \eqref{5.20a} in the above
basis of $\ker P$ and by $n_j$ for $1\le j\le r$ the number of groups of
elements of the form \eqref{5.20b} with $\lambda=\lambda_j$.

Now we give our fundamental definition.

\bde{3}
We say that the wave function $\Psi_W(x,z)$ is a
{\em polynomial Darboux transformation\/} of the Bessel wave function
$\Psi_\beta(x,z)$, $\beta\in\Cset^N$, iff \eqref{3.8} holds (for $V=V_\beta$)
with $P(x,\partial_x)$ and $g(z)$ satisfying:

(i) The kernel of the operator $P$ is homogeneous and
$\Zset_N$-invariant, i.e.\ it has a basis of the form
(\ref{5.20a}, \ref{5.20b}).

(ii) The polynomial $g(z)$ is given by
\beq
g(z)=z^{n_0}\left(z^N-\lambda_1^N\right)^{n_1}\cdots
\left(z^N-\lambda_r^N\right)^{n_r}
\label{5.21}
\eeq
where $n_j$ are the numbers defined above.

We denote the set of all such planes $W$ by $Gr_B(\beta)$ and put
$Gr^{(N)}_B=\bigcup_\beta Gr_B(\beta)$, $\beta\in\Cset^N$-generic.

We say that the polynomial Darboux transformation $\Psi_W(x,z)$ of
$\Psi_\beta(x,z)$ is {\em monomial\/} iff
$$
g(z)=z^{n_0}.
$$
Denote the set of the corresponding planes $W$ by $Gr_{MB}(\beta)$ and put
$Gr^{(N)}_{MB}=\bigcup_\beta Gr_{MB}(\beta)$, $\beta\in\Cset^N$-generic.
\ede

The next theorem provides another equivalent definition of $Gr_B(\beta)$ and
is used essentially in the proof of the bispectrality.
\begin{thm}\label{t5.7}
 The wave function $\Psi_W(x,z)$ is a {\em{polynomial Darboux
transformation\/}} of the Bessel wave function $\Psi_\beta(x,z)$, for generic
$\beta\in\Cset^N$, iff {\em(\ref{3.8}, \ref{3.9}, \ref{5.14},
\ref{5.15})} hold {\em(}for
$V=V_\beta${\em)} and

{\em(i)}  The operator $P$ from {\em \eqref{3.8}}  has the form
\beq
P(x,\partial_x)=x^{-n}\sum_{k=0}^n p_k(x^N)(x\partial_x)^k,
\label{5.22}
\eeq
where $p_k$ are rational functions, $p_n\equiv 1$.

{\em(ii)}  There exists the formal limit
\beq
\lim_{x\to\infty} e^{-xz}\Psi_W(x,z)=1.
\label{5.23}
\eeq
\end{thm}
The limit in \eqref{5.23} is formal in the sense that it is taken in the
coefficients at any power of $z$.

Following G. Wilson \cite{W} we define the {\em{bispectral involution\/}} $b$
on the wave function  $\Psi_W(x,z)$ by exchanging the places of $x$ and $z$:
$$\Psi_{bW}(x,z) = \Psi_W(z,x)$$
(provided the LHS is again a wave function).
\begin{thm}\label{t6.2}
 If $W\in Gr_B(\beta)$ then $bW$ exists and $bW\in
Gr_B(\beta)$.
\end{thm}
This means that $\Psi_{bW}(x,z)$ is a wave function and
\beqa
&&\Psi_{bW}(x,z) =\frac{1}{g_{\rm b}(z)}
P_{\rm b}(x,\partial_x)\Psi_\beta(x,z),
\label{6.4} \\
&&\Psi_\beta(x,z)=\frac{1}{f_{\rm b}(z)} Q_{\rm b}(x,\partial_x)\Psi_{bW}(x,z)
\label{6.9}
\eeqa
for some polynomials $g_{\rm b}, f_{\rm b}$ and operators
$P_{\rm b}, Q_{\rm b}$ satisfying \eqref{5.22}. We can derive explicit
expressions for (\ref{6.4}, \ref{6.9}) in terms of (\ref{3.8}, \ref{3.9}) as
follows.
If the operators $P(x,\partial_x)$ and $Q(x,\partial_x)$ are written in the form
\beqa
&&P(x,\partial_x)=\frac{1}{x^n p_n(x^N)}
\sum_{k=0}^n p_k(x^N)(x\partial_x)^k,
\label{6.5} \hfill\\
&&Q(x,\partial_x)= \sum_{s=0}^m (x\partial_x)^s q_s(x^N)
\frac{1}{x^m q_m(x^N)}
\label{6.11}\hfill
\eeqa
with polynomials $p_k, q_s$ then
\beqa
&&P_{\rm b}(x,\partial_x) = \frac1{g(x)}
\sum_{k=0}^n (x\partial_x)^k p_k\bigl(L_\beta(x,\partial_x)\bigr),
\label{6.6}\hfill\\
&&g_{\rm b}(z) = z^n p_n(z^N)
\label{6.7} \hfill
\eeqa
and
\beqa
&&Q_{\rm b}(x,\partial_x) = \sum_{s=0}^m q_s \left(
L_\beta(x,\partial_x) \right) (x\partial_x)^s \frac{1}{f(x)},
\label{6.13}\hfill\\
&&f_{\rm b}(z) = z^m q_m(z^N).
\label{6.14}\hfill
\eeqa

An immediate corollary is the following result, which we state as a theorem
because of its fundamental character.
\begin{thm}\label{t6.3}
If  $W\in Gr^{(N)}_B$ then the wave function
$\Psi_W(x,z)$ solves the bispectral problem, i.e.\ there exist operators
$L(x,\partial_x)$ and $\Lambda(z,\partial_z)$ such that
\beqa
&& L(x,\partial_x)\Psi_W(x,z) = h(z^N)\Psi_W(x,z),
\label{6.15} \hfill\\
&& \Lambda(z,\partial_z)\Psi_W(x,z) = \Theta(x^N)\Psi_W(x,z).
\label{6.16} \hfill
\eeqa
Moreover,
\beq
\rank A_W=\rank A_{b W}=N.
\label{6.17}
\eeq
\end{thm}
The operators and the polynomials from \eqref{6.15} and \eqref{6.16} are given
by:
\beqa
&&L(x,\partial_x)=P(x,\partial_x)Q(x,\partial_x),
\quad h(z^N)=f(z)g(z);
\label{6.18} \hfill\\
&&\Lambda(z,\partial_z)=P_{\rm b}(z,\partial_z)Q_{\rm b}(z,\partial_z),
\quad \Theta(x^N)=f_{\rm b}(x)g_{\rm b}(x).
\label{6.19} \hfill
\eeqa
The whole bispectral algebra is given in the following proposition.
\begin{prop}\label{p3.5}
For $W\in Gr_B(\beta)$ with generic $\beta\in\Cset^N$ we have
\beq
A_W=\left\{ u\in \Cset[z^N] \mid u(L_\beta)\ker P\subset \ker P\right\}
\label{3.18}
\eeq
and
\beq
\hspace*{-1.85cm}\A_W=\left\{ Pu(L_\beta)P^{-1} \mid u\in A_W\right\}.
\label{3.19}
\eeq
\end{prop}

\bre{1}
$Gr_B^{(1)}$ coincides with the ``adelic'' Grassmannian $Gr^{ad}$ introduced by
Wilson \cite{W}.
The rank 2 bispectral algebras containing an operator of order 2 (the ``even
case'' of Duistermaat and Gr\"unbaum \cite{DG}) are obtained from
$Gr_{MB}^{(2)}\cap Gr^{(2)}$.
\ere

\bre{2}
The (generalized higher) {\em{Airy wave function\/}} is defined by the
following equations (see e.g.\ \cite{Dij})
\beqa
&& L_{\alpha}(x,\partial_x) \Psi_\a(x,z) = z^N \Psi_\a(x,z), \nn \hfill\\
&& L_{\alpha}(x,\partial_x) = \partial_x^N - \alpha_0 x + \sum_{i=2}^{N-1}
\alpha_i \partial_x^{N-i}, \nn \hfill
\eeqa
where $\alpha =
(\alpha_0,\alpha_2,\alpha_3,\ldots,\alpha_{N-1}) \in\Cset^{N-1}$.
The definition of {\em{polynomial Darboux transformations\/}} is similar to
that in the Bessel case (see \thref{5.7})
with only minor modifications: $P$ is not necessarily $\Zset_N$-invariant and
$g(z)$ has to belong to $\Cset[z^N]$. After a suitable definition of the
bispectral involution $b$ we proved analogs of theorems \ref{t6.2} and
\ref{t6.3} (see \cite{BHY2}).
\ere

\bre{3}
The above eigenfunction $\Psi_W(x,z)$ from (\ref{6.15}, \ref{6.16}) is a formal
series. However, if we substitute
$\Phi_\beta(x,z) =\Phi_\beta(xz)$
for $\Psi_\beta(x,z)$ then
$$\Phi_W(x,z)=\frac{1}{g(z)} P(x,\partial_x) \Phi_\beta(x,z)$$
gives a convergent solution to the bispectral problem with the same operators
$L(x,\partial_x)$ and $\Lambda(z,\partial_z)$ (see \cite{BHY2}).
\ere
\section{ $\W$ and monomial Darboux transformations}
The algebra $w_{\infty}$ of the additional symmetries of the KP--hierarchy is
isomorphic to the Lie algebra of regular polynomial differential operators on
the circle:
$$\D = \span \{ z^{\a} \partial_z^{\beta}| \; \a,  \beta
\in \Zset, \; \beta \geq 0 \}. $$
Its unique central extension \cite{KP1, KR} will be denoted by $W_{1+\infty}$.
This algebra gives the action of the additional symmetries on the tau-functions
\cite{ASvM}.

Denote by $c$ the central element of $W_{1+\infty}$ and by $W(A)$
the image of $A \in \D$ under the natural embedding $\D \hookrightarrow
W_{1+\infty}$ (as vector spaces). The algebra $W_{1+\infty}$ has a basis
$$c, \; L_k^l = W(-z^k D^l), \qquad l,k \in \Zset, \; l \geq 0 $$
where $D \equiv D_z = z \partial_z.$ The commutation relations of
$W_{1+\infty}$ can be written most conveniently in terms of generating series
\cite {KR}
$$
\Bigl[ W(z^k e^{xD}),W(z^m e^{yD}) \Bigr] = ( e^{xm} - e^{yk}) W(z^{k+m}
e^{(x+y)D})+ \delta_{k,-m} \frac{e^{xm}-e^{yk}}{1-e^{x+y}}c.
$$

We introduce the subalgebra $\WN$ of $\W$ spanned by $c$ and $L_{kN}^l,$
$l,k\in \Zset,$ $l \geq 0 $. It is a simple fact that $\WN$ is isomorphic to
$\W$.

In the next theorem we sum up some of the results
  from \cite{BHY} which will be needed in Theorems 5 and 6.
\begin{thm}\label{t2.5}
The functions $\tau_\beta(t)$ satisfy the constraints
    \beqa
&& L_0^l \tau_\beta = \lambda_\beta (L_0^l) \tau_\beta, \quad l\geq 0,
\nn\\
&& L_{kN}^l \tau_\beta = 0, \quad k>0, \; l\geq 0,
\nn\\
&& W\Bigl( z^{-kN} P_{\beta,k}(D) D^l \Bigr) \tau_\beta = 0,
\quad k>0, \; l\geq 0,
\nn
    \eeqa
where $P_{\beta,k}(D) = P_\beta(D) P_\beta(D-N) \cdots P_\beta(D-N(k-1)).$
\end{thm}
The first two constraints mean that $\tau_\beta$ is a highest weight vector
with highest weight $\lambda_\beta$ of a representation of $\WN$ in the module
    \beq
\M_\beta = \span \Bigl\{ L_{k_1 N}^{l_1} \cdots L_{k_p N}^{l_p} \tau_\beta
\Big| k_1 \leq\ldots \leq k_p < 0 \Bigr\}.
\label{2.22}
    \eeq
In \cite{BHY} we studied $\M_\beta$ as modules of $\W$. We proved that they are
{\em{quasifinite\/}} (see \cite{KR}) and we derived formulae for the highest
weights and for the singular vectors.

The next theorem establishes the connection between the
highest weight modules $\M_\beta$ and the monomial Darboux transformations.
\begin{thm}\label{t7.1}
If $\tau_W$ is a tau-function lying in the $\WN$-module
$\M_\beta$ {\rm(}$\beta \in \Cset^N${\rm)} then the corresponding plane $W \in
Gr_{MB}(\beta)$. Conversely, if $W \in Gr_{MB}(\beta)$ then $\tau_W \in
\M_{\beta'}$
for some $\beta' \in \Cset^N$ such that $V_{\beta'} \in Gr_{MB}(\beta)$.
\end{thm}
In general $\beta' \not= \beta$. A more precise version of the second part of
\thref{7.1} is given in \cite{BHY3}. Here we shall restrict ourselves only to
the
case when there are no logarithms in the basis \eqref{5.20a} of $\ker P$,
 i.e.\ when it is of the form
\beq
f_k(x) = \sum_{i=1}^{dN} a_{ki} x^{\gamma_i},\quad
0\le k\le n-1,\quad
\label{7.4}
\eeq
where $\gamma=\beta^d$ (see \eqref{5.1'5}).

Let $W\in Gr_{MB}(\beta)$ be a monomial Darboux transformation of a Bessel
plane $V_\beta$, $\beta\in\Cset^N$ with $g(z)=z^n$ and $\ker P$ of the
above type. \deref{3} in this case is equivalent to
\beq
\gamma_i-\gamma_j \in N\Zset \setminus 0
\quad {\rm if}\ a_{ki}a_{kj}\not=0, i\not=j.
\label{7.4'5}
\eeq
We say that the element $f_k(x)$ of the above basis of $\ker P$ is
{\em{associated}} to $\beta_s$ ($1\leq s \leq N$) iff
\beq
\gamma_i-\beta_s\in N\Zset_{\ge0}\quad {\rm if}\; a_{ki}\not=0.
\label{7.7}
\eeq
Then up to a relabeling we can take a subset $\{ \beta_s\}_{1 \leq s\leq M}$
such that
\beq
\beta_s-\beta_t\not\in N\Zset\quad {\rm for}\ \  1\le s\not=t \le M
\label{7.6}
\eeq
and each element of the basis \eqref{7.4} of $\ker P$ is associated to some
$\beta_s$ from this set. Denote by $n_s$ the number of elements associated
to $\beta_s$ and set $n_s=0$ for $s>M$.
Then $n_1+\cdots+n_N = n$. We put
\beq
\beta'=(\beta_1+n_1N-n, \beta_2+n_2N-n,\ldots,\beta_N+n_NN-n).
\label{7.8}
\eeq
\begin{thm}\label{t7.4}
Let $W$ be a monomial Darboux transformation of the Bessel plane $V_\beta$
with $\ker P$ satisfying {\rm(\ref{7.4}, \ref{7.4'5})} and $\beta'$ be as
above.
Then the tau-function  $\tau_W$ of $W$ lies in the $W_{1+\infty}(N)$-module
$\M_{\beta'}$.
\end{thm}
\bex{1}
Let $A=(0\ 1\ 0\ \ldots\ 0)$ and
$\beta'_2-\beta_1'=N\alpha$, $\alpha\in\Zset_{\ge0}$.
Set
$$
\beta''=(\beta_1'-N, \beta_2'+N, \beta_3', \ldots, \beta_N').
$$
Then the module $\M_{\beta''}$ embeds in $\M_{\beta'}$. The singular
vector $\tau_{\beta''}$ is given by
\beq
\tau_{\beta''} = W(P_1)\tau_{\beta'}+\const\cdot\tau_{\beta'},
\label{7.26}
\eeq
where
\beq
P_1 =
-N(\alpha+1) z^{-N}\frac{P_{\beta'}(D_z)} {D_z-\beta_1'}\left(z^{-N}
P_{\beta'} (D_z)\right)^{\alpha+1}.
\label{7.24}
\eeq
\eex

Another important question posed by Duistermaat and Gr\"unbaum \cite{DG} is
about the existence of an hierarchy of symmetries leaving the manifold of
bispectral operators invariant. The following two theorems answer this question
for the manifold of monomial Darboux transformations.
\begin{thm}\label{t7.7}
Vector fields corresponding to $W_{1+\infty}^+(N)$ are
tangent to the manifold $Gr_{MB}^{(N)}$ of monomial Darboux transformations.
More precisely, if $W\in Gr_{MB}^{(N)}$ then
$$
\exp\left(\sum_{i=1}^p \lambda_i L_{Nk_i}^{l_i}\right)\tau_W
$$
is a tau-function associated to a plane from $Gr_{MB}^{(N)}$ for
arbitrary
$p\in \Nset$, $\lambda_i\in\Cset$, $l_i,k_i\in\Zset_{\ge0}$.
\end{thm}
Let us define
$$
\overline L_m
=\frac1N\sum_{k\in{\Zset\setminus N\Zset}} \normprod{ J_{mN-k}J_k }
$$
where $J_k=L^0_k=W(-z^k)$.
The operators $\overline L_m$, $m\in \Zset$ form a Virasoro algebra with
central charge $N-1$ which we denote by $Vir_N$.
Denote by $Vir_N^+$ the subalgebra spanned by $\overline L_m$, $m\ge0$. Then we
can formulate the following theorem which for $N=2$ contains Magri--Zubelli's
result \cite{MZ}.
\begin{thm}\label{t7.8}
The manifold $Gr_B^{(N)}\cap Gr^{(N)}$ is preserved by
the vector fields corresponding to $Vir_N^+$. More precisely, if $W\in
Gr_B^{(N)}\cap Gr^{(N)}$ then
$$
\exp\left(\sum_{i=1}^p \lambda_i\overline L_{k_i}\right) \tau_W
$$
is a tau-function associated to a plane from $Gr_B^{(N)}\cap Gr^{(N)}$ for
arbitrary $p\in \Nset$, $\lambda_i\in\Cset$, $k_i\ge0$.
\end{thm}

\bre{4}
$Gr_B^{(N)}\cap Gr^{(N)}=Gr_{MB}^{(N)}\cap Gr^{(N)}.$
\ere
\section{Explicit formulae and examples}
We shall begin with stating explicit formulae for the bispectral operators from
(\ref{6.15}, \ref{6.16}) in the case of monomial Darboux transformations when
there are
no logarithms in the basis \eqref{5.20a} of $\ker P$. The general case of
monomial Darboux
transformation can be reduced to this one by taking a limit in all formulae
(see \cite{BHY3}).

Let $\beta\in \Cset^N$ and $W\in Gr_{MB}(\beta)$. We use the notation from
(\ref{3.8}, \ref{3.9}, \ref{6.4}, \ref{6.9}).
Then $h(z)=z^d$, $g(z)=z^n$, $f(z)=z^{dN-n}$
for some $n$, $d$. The kernel of $P$ has a basis of the form
(\ref{7.4}, \ref{7.4'5})
where $\gamma=\beta^d$ is from \eqref{5.1'5}.
We shall use multi-index notation for
subsets $I=\{i_1<\ldots< i_n\}$ of $\{1,\ldots,dN\}$.
Let $A$ be the matrix $(a_{ki})$ and $A^I=(a_{k,i_l})_{0\le k,\; l\le n-1}$ be
the corresponding minor of $A$. Denote by $I^0$ the complement of $I$.
We put $\gamma_I = \{\gamma_i\}_{i\in I}$;
$(\delta_I)_i = 1$ for $i\in I$ and $=0$ for $i\in I^0$;
$\Delta_I=\prod_{r<s}(\gamma_{i_r}-\gamma_{i_s})$.
Let $I_{\min}$ be the subset of $\{1,\ldots,dN\}$ with $n$ elements such that
$\det A^{I_{\min}}\not=0$ and $\sum_{i\in I_{\min}} \gamma_i$ be the minimum of
all such sums, and set
$p_I=\sum_{i\in I}\gamma_i -\sum_{i\in I_{\min}}\gamma_i$.
Eq. \eqref{7.4'5}  implies that these numbers are divisible by $N$.
\begin{prop}\label{p9.1}
In the above notation the operators and the polynomials from {\rm{(\ref{6.18},
\ref{6.19})}} are given by the following formulae{\rm{:}}
\begin{eqnarray*}
&&
\hspace*{-1cm}
{\rm{(a)}} \;\;\textstyle  g(z)=z^n,
\hfill\\
&&P= \Bigl( \sum\det A^I\Delta_I x^{p_I} \Bigr)^{-1} \Bigl( \sum\det
A^I\Delta_I x^{p_I}L_{\gamma_I}\Bigr).
\hfill\\
&&
\hfill\\
&&
\hspace*{-1cm}
{\rm{(b)}} \;\;\textstyle  f(z)= z^{dN-n},
\hfill\\
&&Q= \Bigl( \sum\det A^I\Delta_I L_{\gamma_{I^0} -n\delta_{I^0}} x^{p_I} \Bigr)
\Bigl( \sum \det A^I\Delta_I x^{p_I} \Bigr)^{-1}.
\hfill\\
&&
\hfill\\
&&
\hspace*{-1cm}
{\rm{(c)}} \;\;\textstyle  g_{\rm b}(z) =z^n \sum\det A^I\Delta_I
z^{p_I},
\hfill\\
&&P_{\rm b} =\sum\det A^I\Delta_I L_{\gamma_I} (L_\beta)^{p_I/N}.
\hfill\\
&&
\hfill\\
&&
\hspace*{-1cm}
{\rm{(d)}} \;\;\textstyle f_{\rm b}(z)= z^{dN-n} \sum\det A^I\Delta_I
z^{p_I},
\hfill\\
&&Q_{\rm b}=\sum\det A^I\Delta_I (L_\beta)^{p_I/N}
L_{\gamma_{I^0}-n\delta_{I^0}}.
\end{eqnarray*}
\end{prop}

\bex{2}
Consider the case when $\beta^d=\gamma$ has different coordinates. Then choose
the following basis of $\ker L_\beta^d$:
$$
\Phi_{(k-1)d+j}(x) =  \mu_{kj} x^{\beta_k+(j-1)N}, \quad 1\le k\le N,\ 1\le
j\le d,
$$
where
$$\mu_{k,1} := 1, \quad \mu_{kj} :=  \mu_{k,j-1}\cdot\prod_{i=1}^N
(\beta_i-\beta_k-(j-1)N)^{-1}.$$
In this basis the action of $L_\beta$ is quite simple:
$L_\beta\Phi_{(k-1)d+j} = \Phi_{(k-1)d+j-1}$ for $2\le j\le d$ and
$=0$ for $j=1$.
Let a basis of $\ker P$ be
$$
f_k(x)=\sum_{i=1}^{dN} a_{ki}\Phi_i(x),\quad k=0,\ldots,d-1.
$$
Let $\beta_i-\beta_j\in N\Zset$ for all $i$, $j$ and
the matrix $A=(a_{ki})$ has the form:
$$
A=\pmatrix{
t_0^{(1)}    &       &       &          &\ldots  & t_0^{(N)}\cr
\noalign{\vskip3pt}
t_1^{(1)}    & t_0^{(1)} &   &          &\ldots  & t_1^{(N)} &t_0^{(N)}\cr
\noalign{\vskip3pt}
t_2^{(1)}    & t_1^{(1)} & t_0^{(1)} &   &\ldots  & t_2^{(N)} & t_1^{(N)}
   & t_0^{(N)}\cr
\vdots  & \vdots &\hfill&\ddots\hfill & &   \vdots & \vdots &\ddots\hfill\cr
\noalign{\vskip3pt}
t_{n-1}^{(1)} & t_{n-2}^{(1)} &\ldots& t_0^{(1)} & \ldots  &
t_{n-1}^{(N)} & t_{n-2}^{(N)} &\ldots& t_0^{(N)}}
$$
Then the operator $L=P L_\beta P^{-1}$ is differential of order $N$ and solves
the bispectral problem. For a generic $\beta\in \Cset^N$ the spectral algebra
has rank $N$ (i.e.\ it is $\Cset [L]$). For $N=2$ this is the {\em{``even
case''\/}} of J.~J.~Duistermaat and F.~A.~Gr\"unbaum \cite{DG} (see also
\cite{MZ}).
\eex

\bex{3}
All {\it bispectral algebras of rank 1\/} are polynomial
Darboux transformations of the plane $H_+=\span\{z^k\}_{k \geq0}$ (see
\cite{W}). This corresponds to the $N=1$ Bessel with
$$
\beta=(0),\quad L_{(0)}=\partial_x, \quad V_{(0)}=H_+,
\quad \tau_{(0)}(t)=1, \quad
\Psi_{(0)}(x,z)=e^{xz}.
$$
The operator $L$ which solves the bispectral problem is a
Darboux transformation of the operator $h\left(L_{(0)}\right) =h(\partial_x)$
with constant coefficients. The ``adelic Grassmannian'' $Gr^{ad}$, introduced
by Wilson \cite{W}, coincides with $Gr_B((0))$ $(=Gr_B^{(1)})$.
\eex

In the last example we study the simplest polynomial Darboux transformation of
a Bessel operator of order 2.

\bex{4}
For $N=2$, $\beta=(1-\nu,\nu)$ the corresponding Bessel
operator is
$$
L_\beta=x^{-2}(D_x-(1-\nu))(D_x-\nu)= \partial_x^2 +\frac{\nu(1-\nu)}{x^2},
\quad D_x=x\partial_x.
$$
Fix $\lambda\in \Cset\setminus\{0\}$ and set
$
h(z)=(z-\lambda^2)^2.
$
For fixed $\lambda, a\in\Cset\setminus\{0\}$ we take $\ker P$ to have a basis
$$
f_k(x) = \Psi_\beta(x,(-1)^k \lambda) + a D_x \Psi_\beta(x,(-1)^k \lambda),
\qquad k=0,1.
$$
After introducing the operator
$$
 P(a,\lambda,\mu) = \frac{1}{x^2 p_2(x^2)}
\Bigl\{ p_2(x^2) D^2_x + p_1(x^2) D_x + p_0(x^2)\Bigr\},
$$
where
$\mu^2=(a+1-a^2\nu(\nu-1))/a^2\lambda^2$, $p_2(x^2)=x^2-\mu^2$,
$p_1(x^2)=\mu^2-3x^2$, $p_0(x^2) = -\lambda^2 x^4 +(2\lambda^2\mu^2
+(a+1)(2a-1)a^{-2})
x^2 +((a+1)a^{-2} -\lambda^2\mu^2)\mu^2$,
the operators and polynomials from formulae
(\ref{6.18}, \ref{6.19}) are given by:
\begin{eqnarray}
&&P = P(a,\lambda,\mu), \;\;\;\qquad g(z) = z^2 - \lambda^2;
\label{100}
\\
&&Q = P^*(b,\lambda,\mu), \;\;\qquad f(z) = z^2 - \lambda^2;
\label{101}
\\
&&P_{\rm b} = P(a,\mu,\lambda), \,\;\qquad g_{\rm b}(z) = z^2 - \mu^2;
\label{102}
\\
&&Q_{\rm b} = P^*(b,\mu,\lambda), \qquad f_{\rm b}(z) = z^2 - \mu^2;
\label{103}
\end{eqnarray}
where $b=-a/(a+1)$ and ``*'' is the formal adjoint of differential
operators (i.e.\
the unique antiautomorphism such that $\partial_x^* = -\partial_x$, $x^* = x).$
 The spectral algebras
$$
\A_W= P\left(L_\beta-\lambda^2\right)^2 \Cset[L_\beta] P^{-1}
$$
and
$$
\A_{bW} =  P_{\rm b}\left(L_\beta-\mu^2\right)^2 \Cset
[L_\beta] P^{-1}_{\rm b}
$$
consist of operators of orders $4,6,8,10,\ldots$

The above formulae (\ref{100}--\ref{103}) raise some interesting questions
about the properties of the bispectral involution and of the so-called adjoint
involution (see \cite{BHY2}).
\eex

\flushleft{\bf{Acknowledgements}}

\medskip\noindent
We are grateful to the organizers of the conference
in Geometry and Mathematical Physics, Zlatograd 95 for the kind hospitality.
This work was partially supported by Grant MM--402/94 of Bulgarian
Ministry of Education, Science and Technologies.
\begin{small}
    
\end{small}

\begin{thebibliography}{9}
    \bibitem{ASvM}
Adler, M., Shiota, T., van~Moerbeke, P.:
{\em{A Lax representation for the vertex operator and the central extension}}.
Commun. Math. Phys. {\bf 171}, 547--588 (1995).
     \bibitem{BHY}
Bakalov, B., Horozov, E., Yakimov, M.: {\em{Tau-functions as highest weight
vectors for $W_{1+\infty}$ algebra}}. Sofia preprint (1995), hep-th/9510211.
     \bibitem{BHY1}
Bakalov, B., Horozov, E., Yakimov, M.: {\em{B\"acklund--Darboux transformations
in Sato's Grassmannian}}. Sofia preprint (1996), q-alg/9602010.
     \bibitem{BHY2}
Bakalov, B., Horozov, E., Yakimov, M.: {\em{Bispectral algebras of commuting
ordinary differential operators}}, Sofia preprint (1996), q-alg/9602011.
     \bibitem{BHY3}
Bakalov, B., Horozov, E., Yakimov, M.: {\em{Highest weight modules over
$W_{1+\infty}$ algebra and the bispectral problem}}. Sofia preprint (1996),
q-alg/9602012.
     \bibitem{BE}
Bateman, H., Erd\'elyi, A.: {\em{Higher transcendental functions\/}}.
New York: McGraw-Hill, 1953.
     \bibitem{BC}
Burchnall, J.L., Chaundy, T.W.: {\em{Commutative ordinary differential
operators}}.
Proc. Lond. Math. Soc. {\bf 21}, 420--440 (1923);
Proc. Royal Soc. London (A) {\bf 118}, 557--583 (1928);
Proc. Royal Soc. London (A) {\bf 134}, 471--485 (1932).
     \bibitem{Da}
Darboux, G.: {\em{Le{\c c}ons sur la th\'eorie g\'en\'erale des surfaces}}.
2\`eme partie, Paris: Gauthiers--Villars, 1889.
     \bibitem{Dij}
Dijkgraaf, R.: {\em{Intersection theory, integrable hierarchies and
topological field theory}}. Lecture Notes at Cargese Summer School (1991),
hep-th/9201003.
     \bibitem{DG}
Duistermaat, J.J., Gr\"unbaum, F.A.: {\em{Differential equations in the
spectral parameter}}.
Commun. Math. Phys. {\bf 103}, 177--240 (1986).
     \bibitem{F}
Fastr\'e, J.: {\em{B\"acklund--Darboux transformations and $W$-algebras}}.
Doctoral Dissertation, Univ. of Louvain, 1993.
     \bibitem{G1}
Gr\"unbaum, F.A.: {\em{The limited angle reconstruction problem in computer
tomography}}.
Proc. Symp. Appl. Math. {\bf 27}, AMS, L. Shepp (ed.), 43--61 (1982).
      \bibitem{G3}
Gr\"unbaum, F.A.: {\em{Time-band limiting and the bispectral problem}}.
Comm. Pure Appl. Math. {\bf 47}, 307--328 (1994).
     \bibitem{HH}
Haine, L., Horozov, E.: {\em{Tau-functions and modules over the Virasoro
algebra}}. In: {\em{Abelian varieties}}, W. Barth et al. (eds.), Berlin,
New York: Walter de Gruyter, 1995.
     \bibitem{I}
Ince, E.L.: {\em{Ordinary Differential Equations}}. New York: Dover, 1944.
     \bibitem{KP1}
Kac, V.G., Peterson, D.H.: {\em{Spin and wedge representations of
infinite-dimensional Lie algebras and groups}}.
Proc. Natl. Acad. Sci. USA {\bf 78}, 3308--3312 (1981).
      \bibitem{KR}
Kac, V.G., Radul, A.: {\em{Quasifinite highest weight modules over the Lie
algebra of differential operators on the circle}}.
Commun. Math. Phys. {\bf 157}, 429--457 (1993), hep-th/9308153.
    \bibitem{KRa}
Kac V.G., Raina A.: {\em{Bombay lectures on highest
weight representations of infinite dimensional Lie algebras}}. Adv. Ser.
Math. Phys. {\bf2}, Singapore: World Scientific, 1987.
     \bibitem{KrN}
Krichever, I., Novikov, S.: {\em{Holomorphic bundles over algebraic
curves and nonlinear equations}}.
Russian Math. Surveys {\bf 35}, 53--79 (1980).
     \bibitem{MZ}
Magri, F., Zubelli, J.: {\em{Differential equations in the
spectral parameter, Darboux transformations and a hierarchy of
master equations for KdV}}.
Commun. Math. Phys. {\bf 141}, 329--351 (1991).
     \bibitem{S}
Sato, M.: {\em{Soliton equations as dynamical systems on infinite
dimensional Grassmann manifolds}}.
RIMS Kokyuroku {\bf 439}, 30--40 (1981).
     \bibitem{SW}
Segal, G., Wilson, G.: {\em{Loop Groups and equations of KdV type}}.
Publ. Math. IHES {\bf 61}, 5--65 (1985).
    \bibitem{W}
Wilson, G.: {\em{Bispectral commutative ordinary differential operators}}.
J. Reine Angew. Math. {\bf 442}, 177--204 (1993).
    \bibitem{Z}
Zubelli, J.: {\em{Differential equations in the spectral parameter
for matrix differential operators}}.
Physica {\bf D 43}, 269--287 (1990).
    \end{thebibliography}
\end{document}